\documentclass[11pt]{article}
\usepackage{amssymb,latexsym}
\usepackage{epsfig}
\usepackage{eufrak}
\usepackage{amsmath}
\usepackage{mathrsfs}
\usepackage{color}

\newtheorem{theorem}{Theorem}
\newtheorem{lemma}{Lemma}

\newcommand{\be}{\begin{equation}}
\newcommand{\ee}{\end{equation}}
\newcommand{\bee}{\begin{eqnarray*}}
\newcommand{\eee}{\end{eqnarray*}}
\newcommand{\bel}{\begin{eqnarray}}
\newcommand{\eel}{\end{eqnarray}}
\newcommand{\bec}{\begin{cases}}
\newcommand{\eec}{\end{cases}}
\newcommand{\bem}{\begin{bmatrix}}
\newcommand{\eem}{\end{bmatrix}}

\newcommand{\la}{\label}
\newcommand{\li}{\left}
\newcommand{\ri}{\right}

\newcommand{\ovl}{\overline}
\newcommand{\udl}{\underline}

\newcommand{\lf}{\lfloor}
\newcommand{\rf}{\rfloor}

\newcommand{\ep}{\epsilon}

\newcommand{\lm}{\lambda}

\newcommand{\de}{\delta}
\newcommand{\De}{\Delta}

\newcommand{\vse}{\vartheta}
\newcommand{\se}{\theta}
\newcommand{\Se}{\Theta}

\newcommand{\al}{\alpha}
\newcommand{\ba}{\beta}

\newcommand{\Om}{\Omega}

\newcommand{\f}{\frac}

\newcommand{\cd}{\cdots}

\newcommand{\qu}{\quad}
\newcommand{\qqu}{\qquad}
\newcommand{\fa}{\forall}

\newcommand{\mscr}{\mathscr}
\newcommand{\mcal}{\mathcal}

\newcommand{\mrm}{\mathrm}

\newcommand{\ap}{\approx}

\newcommand{\sh}{\slash}

\newcommand{\tx}{\text}

\newcommand{\iy}{\infty}

\newcommand{\bed}{\begin{description}}
\newcommand{\eed}{\end{description}}
\newcommand{\bei}{\begin{itemize}}
\newcommand{\eei}{\end{itemize}}
\newcommand{\ben}{\begin{enumerate}}
\newcommand{\een}{\end{enumerate}}

\newcommand{\beL}{\begin{lemma}}
\newcommand{\eeL}{\end{lemma}}
\newcommand{\beT}{\begin{theorem}}
\newcommand{\eeT}{\end{theorem}}

\newcommand{\bpf}{\begin{pf}}
\newcommand{\epf}{\end{pf}}
\newcommand{\bsk}{\bigskip}
\newcommand{\bi}{\binom}

\newcommand{\pfbox}{\hfill\mbox{$\Box$}}

\newenvironment{pf}{\paragraph*{Proof{\rm.}}}{\pfbox\bigskip}

\begin{document}

\title{{\bf  Coverage Probability of Random Intervals}
\thanks{The author had been previously working with Louisiana State University at Baton Rouge, LA 70803, USA,
and is now with Department of Electrical Engineering, Southern
University and A\&M College, Baton Rouge, LA 70813, USA; Email:
chenxinjia@gmail.com. }}

\author{Xinjia Chen}

\date{July 2007}

\maketitle

\begin{abstract}

In this paper, we develop a general theory on the coverage
probability of random intervals defined in terms of discrete random
variables with continuous parameter spaces. The theory shows that
the minimum coverage probabilities of random intervals with respect
to corresponding parameters are achieved at discrete finite sets and
that the coverage probabilities are continuous and unimodal when
parameters are varying in between interval endpoints. The theory
applies to common important discrete random variables including
binomial variable, Poisson variable, negative binomial variable and
hypergeometrical random variable. The theory can be used to make
relevant statistical inference more rigorous and less conservative.

\end{abstract}

\section{Binomial Random Intervals} \la{sec1}

Let $X$ be a Bernoulli random variable defined in a probability space $(\Om, \mscr{F}, \Pr )$ such that $\Pr \{
X = 1\} = p$ and $\Pr \{X = 0 \} = 1 - p$ where $p \in (0,1)$. Let $X_1, \cd, X_n$ be $n$ identical and
independent samples of $X$. In many applications, it is important to construct a confidence interval $(L, U)$
such that $\Pr \{ L < p < U  \mid p \} \ap 1 - \de$ with $\de \in (0, 1)$.  Here $L = L(n, \de, K)$ and $U =
U(n, \de, K)$ are multivariate functions of $n, \; \de$ and random variable $K = \sum_{i =1}^n X_i$.  To simply
notations, we drop the arguments and write $L = L(K)$ and $U = U(K)$.  Also, we use notation $\Pr \{ L(K) < p <
U(K) \mid p \}$ to represent the probability when the binomial parameter assumes value $p$.   Such notation is
used in a similar way throughout this paper.  We would thus advise the reader to distinguish this notation from
conventional notation of conditional probability.

\bsk

 Clearly, the construction of confidence interval is independent
of the binomial parameter $p$. But, for fixed $n$ and $\de$, the
quantity $\Pr \{ L (K) < p < U (K) \mid p \}$ is a function of $p$
and is conventionally referred to as the coverage probability.   In
many situations, it is desirable to know what is the worst-case
coverage probability for $p$ belonging to interval $[a, b] \subset
(0,1)$.  For this purpose, we have

\beT \la{thm18} Suppose that both $L(k)$ and $U(k)$ are monotone
functions of $k \in \{0, 1, \cd, n \}$. Then, the minimum of $\Pr \{
L (K) < p < U (K) \mid p \}$ with respect to $p \in [a, b]$ is
attained at the discrete set $\{a, b \} \cup \{ L(k) \in (a, b): 0
\leq k \leq n \} \cup \{U(k) \in (a, b): 0 \leq k \leq n\}$. \eeT

We would like to emphasis that the only assumption in Theorem 1 is
that both $L(k)$ and $U(k)$ are either non-decreasing or
non-increasing with respect to $k$. The interval $(L(K), U(K))$ can
be general {\it random interval} without being restricted to the
context of confidence intervals. Theorem \ref{thm18} can be
generalized as Theorem \ref{Fundamental} in Section \ref{sec4}.  The
application of the
 theorem is discussed in the full version of our paper \cite{Chen4}. Specially, Theorem \ref{thm18} can be
 applied to the sample size problems studied in \cite{Chen1}.

For closed confidence interval $[L,U]$, it is interesting to compute
the infimum of $\Pr \{ L (K) \leq p \leq U (K) \mid p \}$ with
respect to $p \in [a, b] \subset (0,1)$.  For this purpose, we have

\beT  \la{thm288} Suppose that both $L(k)$ and $U(k)$ are monotone
functions of $k \in \{0, 1, \cd, n \}$.  Then, the infimum of $\Pr
\{ L (K) \leq p \leq U (K) \mid p \}$ with respect to $p \in [a, b]$
equals the minimum of the set $\{C(a), \; C(b) \} \cup \li \{ C_U
(p): p \in \mscr{S}_U
 \ri \} \cup \li \{  C_L (p): p \in \mscr{S}_L  \ri
\}$, where {\small \[ \mscr{S}_U = \{U(k) \in (a, b): 0 \leq k \leq
n\}, \qu \mscr{S}_L = \{ L(k) \in (a, b): 0 \leq k \leq n \}, \qu
C(p) = \Pr \{ L (K) \leq p \leq U (K) \mid p \},
\]}
{\small $ C_U (p) = \Pr \{ L (K) \leq p < U (K) \mid p \}$} and
{\small $ C_L (p) = \Pr \{ L (K) < p \leq U (K) \mid p \}$}. \eeT

It should be noted that the only assumption in the above theorem is
that both $L(k)$ and $U(k)$ are either non-decreasing or
non-increasing with respect to $k$. The interval $[L(K), U(K)]$ can
be general {\it random interval} without being restricted to the
context of confidence intervals.  Theorem \ref{thm288} can be
considered as a specialized result of Theorem \ref{Fundamental} in
Section \ref{sec4}.

\section{Poisson Random Intervals} \la{sec2}

Let $X$ be a Poisson random variable defined in a probability space
$(\Om, \mscr{F}, \Pr )$ such that
\[
\Pr \{ X = k \} = \f{\lm^k e^{- \lm} } {k!}, \qqu  k = 0,1, 2, \cd
\]
where $\lm > 0$ is called the Poisson parameter. Let $X_1, \cd, X_n$
be $n$ identical and independent samples of $X$. It is a frequent
problem to construct a confidence interval $(L, U)$ such that $\Pr
\{ L < \lm < U \mid \lm\} \ap 1 - \de$ with $\de \in (0, 1)$.  Here
$L = L(n, \de, K)$ and $U = U(n, \de, K)$ are multivariate functions
of $n, \; \de$ and random variable $K = \sum_{i =1}^n X_i$. For
simplicity of notations, we drop the arguments and write $L = L(K)$
and $U = U(K)$.  For fixed $n$ and $\de$,  the coverage probability
$\Pr \{ L (K) < \lm < U (K) \mid \lm \}$ is a function of $\lm$. The
worst-case coverage probability with respect to $\lm$ belonging to
interval $[a, b] \subset (0, \iy)$ can be obtained by the following
theorem.

\beT \la{thm28} Suppose that both $L(k)$ and $U(k)$ are monotone
functions of non-negative integer $k$. Then, the minimum of $\Pr \{
L (K) < \lm < U (K) \mid \lm \}$ with respect to $\lm \in [a, b]$ is
attained at the discrete set $\{a, b \} \cup \{ L(k) \in (a, b): k
\geq  0 \} \cup \{U(k) \in (a, b): k \geq 0 \}$.  \eeT

It should be emphasized  that the interval $(L(K), U(K))$ can be
general {\it random interval} without being restricted to the
context of confidence intervals. The only assumption in the above
theorem is that both $L(k)$ and $U(k)$ are either non-decreasing or
non-increasing with respect to $k$.  Theorem \ref{thm28} can be
generalized as Theorem \ref{Fundamental} in Section \ref{sec4}.  The
application of the theorem is discussed in the full version of our
paper \cite{Chen4} for the sample size problems  studied in
\cite{Chen2}.

For the exact computation of the infimum of coverage probability
$\Pr \{ L (K) \leq \lm \leq U (K) \mid \lm \}$ for the closed
confidence interval $[L, U]$, we have

 \beT \la{thm48}  Suppose that both $L(k)$ and $U(k)$ are monotone functions of
non-negative integer $k$.  Then, the infimum of $\Pr \{ L (K) \leq
\lm \leq U (K) \mid \lm \}$ with respect to $\lm \in [a, b]$ equals
the minimum of the set $\{C(a), \; C(b) \} \cup \li \{ C_U (\lm):
\lm \in \mscr{S}_U \ri \} \cup \li \{ C_L (\lm): \lm \in \mscr{S}_L
\ri \}$ where {\small \[ \mscr{S}_U = \{U(k) \in (a, b): k \geq 0
\}, \qqu \mscr{S}_L = \{ L(k) \in (a, b): k \geq 0 \}, \qqu C(\lm) =
\Pr \{ L (K) \leq \lm \leq U (K) \mid \lm \}, \]
\[
C_U (\lm) = \Pr \{ L (K) \leq \lm < U (K) \mid \lm \}, \qqu C_L
(\lm) = \Pr \{ L (K) < \lm \leq U (K) \mid \lm \}.
\]} \eeT

In Theorem \ref{thm48}, the interval $[L(K), U(K)]$ can be general
{\it random interval} without being restricted to the context of
confidence intervals.  This theorem is a special case of Theorem
\ref{Fundamental} in Section \ref{sec4}.

\section{Negative-Binomial Random Intervals} \la{sec3}

Let $K$ be a negative binomial random variable such that \be
\la{neb}
 \Pr \{ K = k \} = \bi{k + r - 1 }{k} p^r (1 - p)^k, \qqu k
= 0, 1, \cd \ee with parameter $p \in (0,1)$ and $r > 0$.  In the
special case that $r = 1$, a negative binomial random variable
becomes a geometrical random variable.  For the coverage probability
of open random interval $(L(K), U(K))$ for a negative binomial
random variable $K$, we have

\beT Suppose that both $L(k)$ and $U(k)$ are monotone functions of
non-negative integer $k$. Then, the minimum of $\Pr \{ L (K) < p < U
(K) \mid p \}$ with respect to $p \in [a, b] \subset (0,1)$ is
attained at the discrete set $\{a, b \} \cup \{ L(k) \in (a, b): k
\geq 0 \} \cup \{U(k) \in (a, b): k \geq 0 \}$. \eeT

This theorem can be readily obtained by applying Theorem 7 of
Section 4. For the coverage probability of closed random interval
$[L(K), U(K)]$ for a negative binomial random variable $K$, we have

\beT  Suppose that both $L(k)$ and $U(k)$ are monotone functions of
non-negative integer $k$. Then, the infimum of $\Pr \{ L (K) \leq p
\leq U (K) \mid p \}$ with respect to $p \in [a, b] \subset (0,1)$
equals the minimum of the set $\{C(a), \; C(b) \} \cup \li \{ C_U
(p): p \in \mscr{S}_U
 \ri \} \cup \li \{  C_L (p): p \in \mscr{S}_L  \ri
\}$, where {\small \[ \mscr{S}_U = \{U(k) \in (a, b): k \geq 0 \},
\qu \mscr{S}_L = \{ L(k) \in (a, b): k \geq 0 \}, \qu C(p) = \Pr \{
L (K) \leq p \leq U (K) \mid p \},
\]}
{\small $ C_U (p) = \Pr \{ L (K) \leq p < U (K) \mid p \}$} and
{\small $ C_L (p) = \Pr \{ L (K) < p \leq U (K) \mid p \}$}. \eeT

This theorem can be easily deduced from Theorem \ref{Fundamental} of
next section.

\section{Fundamental Theorem of Random Intervals} \la{sec4}

In previous sections, we discuss coverage probability of random
intervals for specific random variables.  Actually, the results can
be generalized to a large class of discrete random variables. In
this direction, we have recently established in \cite{Chen4} the
following fundamental theorem of random intervals.

\beT \la{Fundamental} Let $K$ be an integer-valued random variable
parameterized by $\se \in \Se$.  Let $L(K)$ and $U(K)$ be functions
of random variable $K$. Let $[a, b]$ be an interval contained in
$\Se$. Let $\mscr{S}_L$ denote the intersection of the interval $(a,
b)$ and the support of $L(K)$. Let $\mscr{S}_U$ denote the
intersection of the interval $(a, b)$ and the support of $U(K)$.
Suppose that, for any $\vartheta \in \Se$, $\Pr \{ L(K) \leq
\vartheta \leq U(K) \mid \se \}$ is a continuous and unimodal
function of $\se \in \Se$. Then, the minimum of $\Pr \{ L (K) < \se
< U (K) \mid \se \}$ with respect to $\se \in [a, b]$ is attained at
the set $\mscr{S}_L \cup \mscr{S}_U \cup \{a, b \} $ and the infimum
of $\Pr \{ L (K) \leq \se \leq U (K) \mid \se \}$ with respect to
$\se \in [a, b]$ is equal to the minimum of the set $\li \{ C_L
(\se):  \se \in \mscr{S}_L \ri \} \cup \li \{ C_U (\se): \se \in
\mscr{S}_U \ri \} \cup \{C(a), \; C_U(a), \;  C(b), \; C_L (b) \} $,
where $ C_L (\se) = \Pr \{ L (K) < \se \leq U (K) \mid \se \}, \;
C_U (\se) = \Pr \{ L (K) \leq \se < U (K) \mid \se \}$ and $C(\se) =
\Pr \{ L (K) \leq \se \leq U (K) \mid \se \}$. Moreover, for both
open random interval $((L(K), U(K))$ and closed random interval
$[L(K), U(K)]$, the coverage probability is continuous and unimodal
for $\se \in (\se^\prime, \se^{\prime \prime})$, where $\se^\prime$
and $\se^{\prime \prime}$ are arbitrary consecutive distinct
elements of $\mscr{S}_L \cup \mscr{S}_U \cup \{a, b \}$.

\eeT

Theorem \ref{Fundamental} is proved in Appendices \ref{seca}.   The
concepts of support and unimodal functions have been used in Theorem
\ref{Fundamental}.  The support of a random variable is actually the
set of all possible values assumed by that random variable.   A
function is said to be a unimodal function of $\se \in \Se$ if there
exists $\se^*$ such that the function is non-decreasing for $\se \in
\Se$ no greater than $\se^*$ and non-increasing for $\se \in \Se$ no
less than $\se^*$.  It should be noted that a monotone function can
be considered as a special case of unimodal function by specifying
$\se^*$ as  the infimum or supremum of $\Se$.  Based on such notion
of unimodal function, the coverage theory stated in Theorem
\ref{Fundamental} applies to one-sided random intervals such as $(-
\iy, U(K)], \; [L(K), \iy), \; (- \iy, U(K)), \; (L(K), \iy)$.

Under the assumption that $\{L (K) \leq \vse \leq U(K) \}$ is an
event that $K$ is contained in an interval, it can be readily shown
that the assumption of Theorem \ref{Fundamental} is satisfied for
common discrete random variables such as binomial random variable,
Poisson random variable, geometrical random variable, negative
binomial random variable, etc.

Let $C_L(\se) $ and $C_U(\se)$ be defined as in Theorem
\ref{Fundamental}.  By the same argument as that for proving Theorem
\ref{Fundamental}, we can show that the infimum of $\Pr \{ L (K) <
\se \leq U (K) \mid \se \}$ with respect to $\se \in [a, b]$ is
equal to the minimum of the set $\li \{ C_L (\se): \se \in
\mscr{S}_L \ri \} \cup \li \{ C_U (\se): \se \in \mscr{S}_U \ri \}
\cup \{C(a), \; C_U(a), \; C(b), \; C_L (b) \} $, where $C(\se) =
\Pr \{ L (K) < \se \leq U (K) \mid \se \}$.  Similarly, the infimum
of $\Pr \{ L (K) \leq \se < U (K) \mid \se \}$ with respect to $\se
\in [a, b]$ is equal to the minimum of the set $\li \{ C_L (\se):
\se \in \mscr{S}_L \ri \} \cup \li \{ C_U (\se): \se \in \mscr{S}_U
\ri \} \cup \{C(a), \; C_U(a), \; C(b), \; C_L (b) \} $, where
$C(\se) = \Pr \{ L (K) \leq \se < U (K) \mid \se \}$.

\section{Infimum Coverage Probability over Parameter Space}
\la{infsec}

In previous sections, we have considered the infimum of coverage
probability over a closed interval $[a, b]$ contained in the
parameter space $\Se$. In many cases, the parameter space $\Se$ is
an open set and consequently,  the infimum of coverage probability
over $\Se$ needs to be treated differently.

As an application of Theorem \ref{Fundamental}, we have obtained the
following results for binomial random intervals.

\beT  \la{infcov} Let $K = \sum_{i=1}^n X_i$, where $X_1, \cd, X_n$
are i.i.d. samples of Bernoulli random variable $X$ such that $\Pr
\{ X = 1 \} = 1 - \Pr \{ X = 0 \} = p \in (0, 1)$. Let $L(k)$ and
$U(k)$ be functions of  nonnegative integer $k$ such that $0 = L(0)
< U(0) < 1, \; 0 < L(n) < U(n) = 1$ and that, for any $\se \in (0,
1)$, there exist two numbers $u$ and $v$ such that $\{ L(K) \leq \se
\leq U(K) \} = \{ u \leq K \leq v \}$. Let $\mscr{S}_L = \{ L(k) \in
(0, 1): k = 1, \cd, n  \}, \; \mscr{S}_U = \{U(k)  \in (0, 1): k =
0, 1, \cd, n - 1 \}$ and $\mscr{S} = \mscr{S}_L \cup \mscr{S}_U$.
Then, $\inf_{p \in (0, 1)} \Pr \{ L(K) < p < U(K) \mid p \}$ is
equal to $\min_{p \in \mscr{S}} \Pr \{ L(K) < p < U(K) \mid p \}$.
Moreover, $\inf_{p \in (0, 1)} \Pr \{ L(K) \leq p \leq U(K) \mid p
\}$ is equal to the minimum of $\min_{p \in \mscr{S}_L} \Pr \{ L(K)
< p \leq U(K) \mid p \}$ and $\min_{p \in \mscr{S}_U} \Pr \{ L(K)
\leq p < U(K) \mid p \}$. Furthermore, $\inf_{p \in (0, 1)} \Pr \{
L(K) < p < U(K) \mid p \} = \inf_{p \in (0, 1)} \Pr \{ L(K) \leq p
\leq U(K) \mid p \}$ under additional assumption that $\mscr{S_L}
\cap \mscr{S}_L = \emptyset$. \eeT

See Appendix \ref{infcov_app} for a proof.

Theorem \ref{infcov} reveals a counterintuitive fact. That is, the
infimum of the coverage probability of an open random interval is
not necessarily equals to the infimum of the coverage probability of
the corresponding closed random interval.  This discovery can be
confirmed by investigating random intervals with
\[
L(K) = \max \li \{  \f{K}{n} - \f{1}{n}, 0 \ri \}, \qqu  U(K) = \min
\li \{  \f{K}{n} + \f{1}{n}, 1 \ri \},
\]
where $K$ is defined in  Theorem \ref{infcov}.  For $n = 3$, we can
show by direct computation that
\[
\Pr \{ L(K) < p < U(K) \mid p \} = \bec (1 - p)^3 + 3 p (1 - p)^2 &
\tx{for} \; 0 < p < \f{1}{3},\\
\f{4}{9} & \tx{for} \; p = \f{1}{3},\\
3 p (1 - p) & \tx{for} \; \f{1}{3} < p < \f{2}{3},\\
\f{4}{9} & \tx{for} \; p = \f{2}{3},\\
3 p^2 (1 - p) + p^3 & \tx{for} \; \f{2}{3} < p < 1 \eec
\]
\[
\Pr \{ L(K) \leq p \leq U(K) \mid p \} = \bec (1 - p)^3 + 3 p (1 -
p)^2 &
\tx{for} \; 0 < p < \f{1}{3},\\
\f{26}{27} & \tx{for} \; p = \f{1}{3},\\
3 p (1 - p) & \tx{for} \; \f{1}{3} < p < \f{2}{3},\\
\f{26}{27} & \tx{for} \; p = \f{2}{3},\\
3 p^2 (1 - p) + p^3 & \tx{for} \; \f{2}{3} < p < 1 \eec
\]
and that \bel &  & \inf_{p \in (0, 1)} \Pr \{ L(K) < p < U(K) \mid p
\} = \min_{p \in (0, 1)} \Pr \{ L(K) < p < U(K) \mid p
\} = \f{4}{9}, \nonumber\\
&  & \inf_{p \in (\f{1}{3}, \f{2}{3})} \Pr \{ L(K) < p < U(K) \mid p
\} = \f{2}{3} > \min_{p \in \{  \f{1}{3}, \f{2}{3}  \} }  \Pr \li \{
L(K) < p < U(K) \mid p \ri \} =  \f{4}{9}, \nonumber\\
&  & \inf_{p \in (\f{1}{3}, \f{2}{3})} \Pr \{ L(K) \leq p \leq U(K)
\mid p \} = \f{2}{3} < \min_{p \in \{  \f{1}{3}, \f{2}{3}  \} }  \Pr
\li \{ L(K) \leq p \leq U(K) \mid p \ri \} =  \f{26}{27}, \nonumber\\
& & \inf_{p \in (0, 1)} \Pr \{ L(K) \leq p \leq U(K) \mid p \} =
\f{2}{3} > \inf_{p \in (0, 1)} \Pr \{ L(K) < p < U(K) \mid p \} =
\f{4}{9}. \la{surprise} \eel In particular, (\ref{surprise}) shows
that the infimum of coverage probabilities for the open and closed
random intervals are not equal.  This is quite surprising. The
coverage probabilities $\Pr \{ L(K) < p < U(K) \mid p \}$ and $\Pr
\{ L(K) \leq p \leq U(K) \mid p \}$ are shown by Figure 1 and Figure
2 respectively.

\begin{figure}
\centering
\includegraphics[height=9cm]{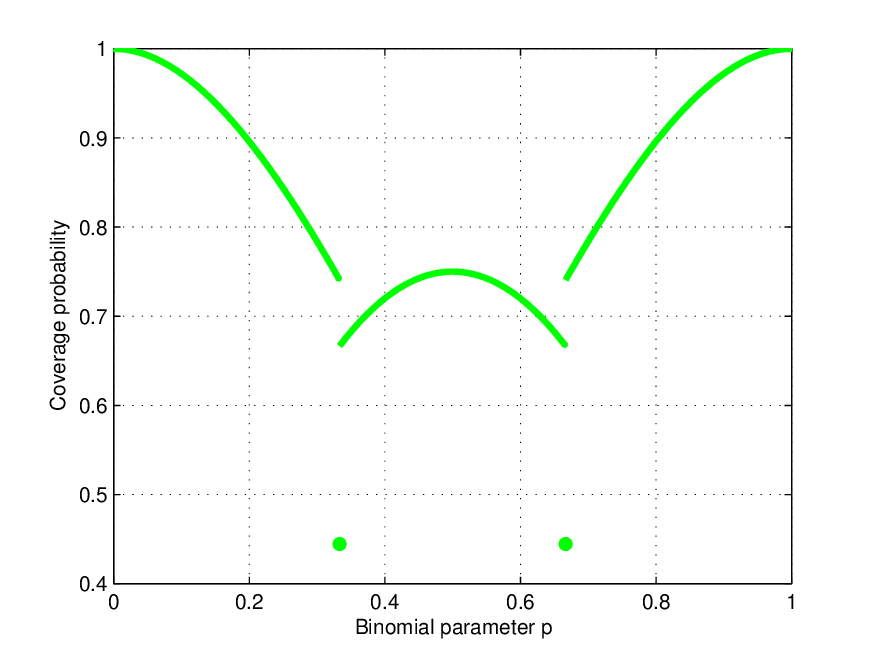}
\caption{Coverage probability of open random interval}
\label{density}       
\end{figure}

\begin{figure}
\centering
\includegraphics[height=9cm]{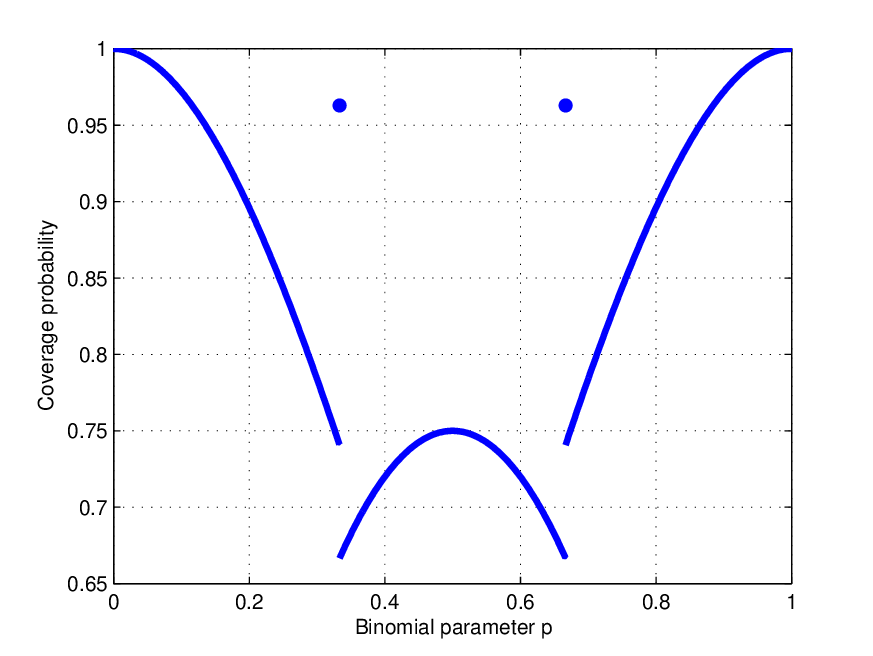}
\caption{Coverage probability of closed random interval}
\label{density}       
\end{figure}

\bsk

The following result establishes the {\it nonexistence} of local
minima for the coverage probability of binomial random intervals
under mild conditions.

\beT \la{nonminima} Let $K = \sum_{i=1}^n X_i$, where $X_1, \cd,
X_n$ are i.i.d. samples of Bernoulli random variable $X$ such that
$\Pr \{ X = 1 \} = 1 - \Pr \{ X = 0 \} = p \in (0, 1)$. Let $L(k)$
and $U(k)$ be nondecreasing functions of nonnegative integer $k$
such that $L(0) = 0, \;  U(n) = 1$ and $L(k) \leq U(k)$ for $k = 0,
1, \cd, n$. Then, there exists no local minima for $\Pr \{ L(K) \leq
p \leq U(K) \mid p \}$ with respect to $p \in (0,1)$. \eeT

The proof of Theorem \ref{nonminima} is available in Appendix
\ref{nonminima_app}.

By similar argument as that for proving Theorems \ref{Fundamental}
and \ref{infcov}, we have established Theorems
\ref{One_BRI}--\ref{neb_infcov} in the sequel.

For one-sided binomial random intervals, we have the following
results.

\beT \la{One_BRI} Let $K = \sum_{i=1}^n X_i$, where $X_1, \cd, X_n$
are i.i.d. samples of Bernoulli random variable $X$ such that $\Pr
\{ X = 1 \} = 1 - \Pr \{ X = 0 \} = p \in (0, 1)$.  Let $L(k)$ and
$U(k)$ be nondecreasing functions of nonnegative integer $k$ such
that $0 = L(0) < L(n) < 1$ and $0 < U(0) < U(n) = 1$. Let
$\mscr{S}_L = \{ L(k) \in (0, 1): k = 1, \cd, n \}$ and $\mscr{S}_U
= \{U(k) \in (0, 1): k = 0, 1, \cd, n - 1 \}$. Then, \bee & &
\inf_{p \in (0, 1)} \Pr \{ L(K) < p \mid p \} = \inf_{p \in (0, 1)}
\Pr \{ L(K) \leq p \mid p \} = \min_{p \in \mscr{S}_L}
\Pr \{ L(K) < p \mid p \},\\
&  & \inf_{p \in (0, 1)} \Pr \{ p < U(K) \mid p \} = \inf_{p \in (0,
1)} \Pr \{ p \leq U(K) \mid p \} = \min_{p \in \mscr{S}_U} \Pr \{ p
< U(K) \mid p \}. \eee
 \eeT

For Poisson random intervals, we have the following results.

 \beT
\la{Poisson_infcov} Let $K$ be a Poisson random variable of mean
$\lm > 0$. Let $L(k)$ and $U(k)$ be functions of nonnegative integer
$k$.  Let $\mscr{S}_L = \{ L(k) \in (0, \iy): k \geq 1 \}, \;
\mscr{S}_U = \{U(k)  \in (0, \iy): k \geq 0 \}$ and $\mscr{S} =
\mscr{S}_L \cup \mscr{S}_U$. Suppose that $0 = L(0) < U(0), \;
\mscr{S}_L \neq \emptyset, \; \sup \mscr{S}_U = \iy$ and that,  for
any $\se \in (0, \iy)$, there exist two numbers $u$ and $v$ such
that $\{ L(K) \leq \se \leq U(K) \} = \{ u \leq K \leq v \}$.
 Then, $\inf_{\lm \in (0, \iy)} \Pr \{ L(K) < \lm < U(K)
\mid \lm \}$ is equal to $\inf_{\lm \in \mscr{S}} \Pr \{ L(K) < \lm
< U(K) \mid \lm \}$. Moreover, $\inf_{\lm \in (0, \iy)} \Pr \{ L(K)
\leq \lm \leq U(K) \mid \lm \}$ is equal to the minimum of
$\inf_{\lm \in \mscr{S}_L} \Pr \{ L(K) < \lm \leq U(K) \mid \lm \}$
and $\inf_{\lm \in \mscr{S}_U} \Pr \{ L(K) \leq \lm < U(K) \mid \lm
\}$.  \eeT

For one-sided Poisson random intervals, we have the following
results.

 \beT
\la{One_Poisson_infcov} Let $K$ be a Poisson random variable of mean
$\lm > 0$.  Let $L(k)$ and $U(k)$ be nondecreasing functions of
nonnegative integer $k$.  Let $\mscr{S}_L = \{ L(k) \in (0, \iy): k
\geq 1 \}$ and $\mscr{S}_U = \{U(k)  \in (0, \iy): k \geq 0 \}$.
Suppose that $0 = L(0) < U(0), \; \mscr{S}_L \neq \emptyset$ and
$\sup \mscr{S}_U = \iy$. Then, \bee & & \inf_{\lm \in (0, \iy)} \Pr
\{ L(K) < \lm \mid \lm \} = \inf_{\lm \in (0, \iy)} \Pr \{ L(K) \leq
\lm \mid
\lm \} = \inf_{\lm \in \mscr{S}_L} \Pr \{ L(K) < \lm \mid \lm \},\\
&  & \inf_{\lm \in (0, \iy)} \Pr \{ \lm < U(K) \mid \lm \} =
\inf_{\lm \in (0, \iy)} \Pr \{ \lm \leq U(K) \mid \lm \} = \inf_{\lm
\in \mscr{S}_U} \Pr \{ \lm < U(K) \mid \lm \}. \eee \eeT

For negative binomial random intervals, we have the following
results.

 \beT
\la{neb_infcov} Let $K$ be a negative binomial random variable
defined by (\ref{neb}).  Let $L(k)$ and $U(k)$ be non-increasing
functions of nonnegative integer $k$.  Let $\mscr{S}_L = \{ L(k) \in
(0, 1): k \geq 1 \}, \; \mscr{S}_U = \{U(k)  \in (0, 1): k \geq 0
\}$ and $\mscr{S} = \mscr{S}_L \cup \mscr{S}_U$. Suppose that $0 <
L(0) < U(0) = 1$ and $\lim_{k \to \iy} L(k) = 0 < \lim_{k \to \iy}
U(k) < 1$.  Then, $\inf_{p \in (0, 1)} \Pr \{ L(K) < p < U(K) \mid p
\}$ is equal to $\inf_{p \in \mscr{S}} \Pr \{ L(K) < p < U(K) \mid p
\}$. Moreover, $\inf_{p \in (0, 1)} \Pr \{ L(K) \leq p \leq U(K)
\mid p \}$ is equal to the minimum of $\inf_{p \in \mscr{S}_L} \Pr
\{ L(K) < p \leq U(K) \mid p \}$ and $\inf_{p \in \mscr{S}_U} \Pr \{
L(K) \leq p < U(K) \mid p \}$.  Furthermore, \bee & & \inf_{p \in
(0, 1)} \Pr \{ L(K) < p \mid p \} = \inf_{p \in (0, 1)} \Pr \{ L(K)
\leq p \mid
p \} = \inf_{p \in \mscr{S}_L} \Pr \{ L(K) < p \mid p \},\\
&  & \inf_{p \in (0, 1)} \Pr \{ p < U(K) \mid p \} = \inf_{p \in (0,
1)} \Pr \{ p \leq U(K) \mid p \} = \inf_{p \in \mscr{S}_U} \Pr \{ p
< U(K) \mid p \}. \eee \eeT

\section{Hypergeometrical Random Intervals} \la{sec5}

So far what we have addressed are random intervals of variables with continuous parameter spaces.  In this
section, we shall consider random intervals when the parameter space is discrete. We focus on the important
hypergeometrical random variable.

 Consider a finite population of $N$ units, among which $M$ units have a
 certain attribute.
 Let $K$ be the number of units found to have the attribute
 in a sample of $n$ units obtained by sampling without replacement.
 The number $K$ is known to be a random variable of hypergeometrical
 distribution.

 It is a basic problem to construct a
 confidence interval $(L, U)$ with $L = L(N,n, \de, K)$ and $U = U(N, n,\de,
 K)$ such that $\Pr \{ L < M < U  \mid M \} \ap 1 - \de$.
 Here, $U$ and $L$ {\it only assume integer values}.  For
 notational simplicity,  we write $L = L(K)$ and $U = U(K)$.
In practice, it is useful to know the minimum of coverage probability $\Pr \{ L < M < U \mid M \}$ with respect
to $M \in [a, b]$, where $a$ and $b$ are integers taken values in between $0$ and $N$.  For this purpose, we
have

\beT \la{fundamental_thm9} Suppose  that $L(0) \leq L(1) \leq \cd
\leq L(n)$ and $U(0) \leq U(1) \leq \cd \leq U(n)$.  Then,  the
minimum of $\Pr \{ L (K) < M < U (K) \mid M \}$ with respect to $M
\in [a, b]$ is attained at the discrete set $I_{UL}$, where $I_{UL}
= \{a, b \} \cup \{ L(k) \in (a, b): 0 \leq k \leq n \} \cup \{U(k)
\in (a, b): 0 \leq k \leq n\}$.  Moreover, $\Pr \{ L (K) < M < U (K)
\mid M \}$ is unimodal  with respect to $M$ in between consecutive
distinct elements of $I_{UL}$.  \eeT

For a proof, see Appendix \ref{secB}.  In Theorem
\ref{fundamental_thm9}, the interval $(L(K), U(K))$ can be general
{\it random interval} without being restricted to the context of
confidence intervals. This theorem can be applied to the sample size
problems discussed in \cite{Chen3}.

\appendix

\section{Proof of Theorem \ref{Fundamental}} \la{seca}

 We need some preliminary results.

 \beL
 \la{simp}
 Suppose that $\{  \se^\prime < L(K) <
\se^{\prime \prime} \} = \{  \se^\prime < U(K) < \se^{\prime \prime}
\} = \emptyset$. Then,
 \[
\{ L(K) < \se < U(K) \} = \{ L(K) \leq \se \leq U(K) \} = \{ L(K)
\leq \se^\prime < U(K) \} = \{ L(K) < \se^{\prime \prime} \leq U(K)
\}
\]
for any $\se \in (\se^\prime, \se^{\prime \prime})$. \eeL

\bpf

By the assumption of the lemma, we have $\{ L(K) < \se^{\prime
\prime} \} = \{ L(K) \leq \se^\prime \} \cup \{ \se^\prime < L(K) <
\se^{\prime \prime} \} = \{ L(K) \leq \se^\prime \}$ and $\{
\se^\prime < L(K) < \se \} \subseteq \{  \se^\prime < L(K) \leq \se
\} \subseteq \{ \se^\prime < L(K) < \se^{\prime \prime} \} =
\emptyset$ for any $\se \in (\se^\prime, \se^{\prime \prime})$.
Consequently, \bel &  & \{ L(K) \leq \se \} = \{ L(K) \leq
\se^\prime \} \cup \{  \se^\prime < L(K) \leq \se \} = \{ L(K) \leq
\se^\prime \} = \{ L(K) <
\se^{\prime \prime} \}, \la{ida}\\
&  & \{ L(K) < \se \} = \{ L(K) \leq \se^\prime \} \cup \{
\se^\prime < L(K) < \se \} = \{ L(K) \leq \se^\prime \} \la{idb}
\eel for any $\se \in (\se^\prime, \se^{\prime \prime})$.  Combining
(\ref{ida}) and (\ref{idb}) yields \be \la{impa}
 \{ L(K) \leq \se
\} = \{ L(K) < \se \} = \{ L(K) \leq \se^\prime \} = \{ L(K) <
\se^{\prime \prime} \}, \qqu \fa \se \in (\se^\prime, \se^{\prime
\prime}). \ee

Again by the assumption of the lemma, we have $\{ U(K) > \se^\prime
\} = \{ U(K) \geq \se^{\prime \prime} \} \cup \{ \se^\prime < U(K) <
\se^{\prime \prime} \} = \{ U(K) \geq \se^{\prime \prime} \}$ and
$\{ \se < U(K) < \se^{\prime \prime} \} \subseteq \{  \se \leq U(K)
< \se^{\prime \prime} \} \subseteq \{ \se^\prime < U(K) <
\se^{\prime \prime} \} = \emptyset$ for any $\se \in (\se^\prime,
\se^{\prime \prime})$. Consequently, \bel &  & \{ U(K) \geq \se \} =
\{ U(K) \geq \se^{\prime \prime} \} \cup \{ \se \leq U(K) <
\se^{\prime \prime}
\} = \{ U(K) \geq \se^{\prime \prime} \} = \{ U(K) > \se^\prime \}, \la{idaa}\\
&  & \{ U(K) > \se \} = \{ U(K) \geq \se^{\prime \prime} \} \cup \{
\se < U(K) < \se^{\prime \prime} \} = \{ U(K) \geq \se^{\prime
\prime} \}  \la{idbb} \eel for any $\se \in (\se^\prime, \se^{\prime
\prime})$. Combining (\ref{idaa}) and (\ref{idbb}) yields \be
\la{impb}
 \{
U(K) \geq \se \} = \{ U(K) > \se \} = \{ U(K) > \se^\prime \} = \{
U(K) \geq \se^{\prime \prime} \}, \qqu \fa \se \in (\se^\prime,
\se^{\prime \prime}). \ee  Taking intersection of events and making
use of (\ref{impa}) and (\ref{impb}), we have \[
 \{ L(K) < \se < U(K) \}
= \{ L(K) \leq \se \leq U(K) \} = \{ L(K) \leq \se^\prime < U(K) \}
= \{ L(K) < \se^{\prime \prime} \leq U(K) \}
\]
for any $\se \in (\se^\prime, \se^{\prime \prime})$.  This completes
the proof of the lemma. \epf

\bsk

Now we are in a position to prove Theorem \ref{Fundamental}.  First,
we shall show the first statement regarding the minimum of $\Pr \{
L(K) < \se < U(K) \mid \se \}$ for $\se \in [a, b]$.  Let
$\se^\prime < \se^{\prime \prime}$ be two consecutive distinct
elements of $\{a, b \} \cup \mscr{S}_U \cup \mscr{S}_L$.   Let $\vse
= \f{ \se^\prime + \se^{\prime \prime}}{2}$.  Then, $\{ \se^\prime <
L(K) < \se^{\prime \prime} \} = \{ \se^\prime < U(K) < \se^{\prime
\prime} \} = \emptyset$ and by Lemma \ref{simp}, we have \bel
 &  & \{ L(K) < \se < U(K)  \} = \{ L(K) \leq \se \leq U(K)  \} \la{imp8}\\
 &  & =
\{ L(K) \leq \vse \leq U(K)  \} = \{ L(K) \leq \se^\prime < U(K) \}
= \{ L(K) < \se^{\prime \prime} \leq U(K) \} \la{imp9} \eel for any
$\se \in (\se^\prime, \se^{\prime \prime})$.  By the assumption of
the theorem,  $\Pr \{ L(K) \leq \vse \leq U(K) \mid \se \}$ is a
continuous and unimodal function of $\se \in \Se$.  It follows from
(\ref{imp8}) and (\ref{imp9}) that both $\Pr \{ L(K) \leq \se^\prime
< U(K) \mid \se \}$ and $\Pr \{ L(K) < \se^{\prime \prime} \leq U(K)
\mid \se \}$ are continuous and unimodal functions of $\se \in
(\se^\prime, \se^{\prime \prime})$. Hence, for
 $\se \in (\se^\prime, \se^{\prime \prime})$,  letting $0 < \ep <
\min (  \se - \se^\prime, \; \se^{\prime \prime} - \se, \; \f{
\se^{\prime \prime} - \se^\prime } {2} )$, we have $\se^\prime + \ep
< \se < \se^{\prime \prime} - \ep$ and {\small \be \la{win1}
 \Pr \{
L(K) < \se^{\prime \prime} \leq U(K) \mid \se \} \geq  \min \li (
\Pr \{ L(K) < \se^{\prime \prime} \leq U(K) \mid \se^\prime + \ep
\}, \; \Pr \{ L(K) < \se^{\prime \prime} \leq U(K) \mid \se^{\prime
\prime} - \ep \} \ri ). \ee}  By virtue of the continuity of $\Pr \{
L(K) < \se^{\prime \prime} \leq U(K) \mid \se \}$ with respect to
$\se \in (\se^\prime, \se^{\prime \prime})$, we have \bel & &
\lim_{\ep \downarrow 0} \Pr \{ L(K) < \se^{\prime \prime} \leq U(K)
\mid \se^\prime + \ep \} = \Pr \{
L(K) < \se^{\prime \prime} \leq U(K) \mid \se^\prime  \}, \la{win2}\\
&  & \lim_{\ep \downarrow 0} \Pr \{ L(K) < \se^{\prime \prime} \leq
U(K) \mid \se^{\prime \prime} - \ep \} = \Pr \{ L(K) < \se^{\prime
\prime} \leq U(K) \mid \se^{\prime \prime} \} \la{win3}\eel It
follows from (\ref{win1}), (\ref{win2}) and (\ref{win3}) that
{\small \be \la{imp18}
 \Pr \{ L(K) < \se^{\prime \prime} \leq U(K) \mid \se
\}\geq \min \li ( \Pr \{ L(K) < \se^{\prime \prime} \leq U(K) \mid
\se^\prime \}, \; \Pr \{ L(K) < \se^{\prime \prime} \leq U(K) \mid
\se^{\prime \prime} \} \ri ) \qqu \ee}  for any $\se \in
(\se^\prime, \se^{\prime \prime})$.  Combining (\ref{imp8}),
(\ref{imp9}) and (\ref{imp18}) yields \bel &  & \Pr \{ L(K) < \se <
U(K) \mid \se \} =  \Pr \{ L(K) \leq \se \leq U(K) \mid \se \}
 = \Pr \{ L(K) < \se^{\prime \prime} \leq U(K) \mid \se \} \qqu \qqu \la{good1}\\
&  & \geq  \min \li (  \Pr \{ L(K) < \se^{\prime \prime} \leq U(K)
\mid \se^\prime \}, \; \Pr \{ L(K) < \se^{\prime \prime} \leq U(K)
\mid \se^{\prime \prime}  \} \ri ) \la{good2}\\
&  &  = \min \li (  \Pr \{ L(K) \leq \se^\prime < U(K) \mid
\se^\prime \}, \; \Pr \{ L(K) < \se^{\prime \prime} \leq U(K) \mid
\se^{\prime \prime}  \} \ri ) \la{good3}\\
&  & \geq   \min \li (  \Pr \{ L(K) < \se^\prime < U(K) \mid
\se^\prime \}, \; \Pr \{ L(K) < \se^{\prime \prime} < U(K) \mid
\se^{\prime \prime}  \} \ri ) \la{vict} \eel for any
 $\se \in (\se^\prime, \se^{\prime \prime})$.
 It can be seen from (\ref{good1}), (\ref{good2}), (\ref{good3}) and (\ref{vict}) that the
 minimum of $\Pr \{ L(K) < \se < U(K) \mid \se \}$ with respect to
 $\se \in [\se^\prime, \se^{\prime \prime}]$ is achieved at either $\se^\prime$ or $\se^{\prime
 \prime}$.    This implies that the minimum of $\Pr \{ L(K) < \se <
U(K) \mid \se \}$ for $\se \in [a, b]$ is attained at $\mscr{S}_L
\cup \mscr{S}_U \cup \{a, b \}$.

Next, we shall show the second statement regarding the infimum of
$\Pr \{ L(K) \leq \se \leq U(K) \mid \se \}$ for $\se \in [a, b]$.
As before, let $\se^\prime < \se^{\prime \prime}$ be two consecutive
distinct elements of $\{a, b \} \cup \mscr{S}_U \cup \mscr{S}_L$.
For simplicity of notations, let \bee &  & \al = \inf_{\se \in
[\se^\prime, \se^{\prime \prime}] } \Pr \{ L(K) \leq \se \leq U(K)
\mid \se \},\\
&  & \ba = \min \li (  \Pr \{ L(K) \leq \se^\prime < U(K) \mid
\se^\prime \}, \; \Pr \{ L(K) < \se^{\prime \prime} \leq U(K) \mid
\se^{\prime \prime}  \} \ri ). \eee  Making use of (\ref{good1}),
(\ref{good2}), (\ref{good3}) and the observation that \[ \ba \leq
\min ( \Pr \{ L(K) \leq \se^\prime \leq U(K) \mid \se^\prime \}, \;
\Pr \{ L(K) \leq \se^{\prime \prime} \leq U(K) \mid \se^{\prime
\prime} \} ), \] we have $\al \geq \ba$.  Now we need to show that
$\al$ is actually equal to $\ba$.  Suppose, to get a contradiction,
that $\alpha$ is greater than $\beta$. Then, \be \la{how}
 \Pr \{ L(K) \leq \se \leq U(K) \mid \se \} > \f{\alpha +
\beta}{2}, \qqu \fa \se \in [ \se^\prime, \se^{\prime \prime} ]. \ee
As a consequence of (\ref{imp8}), (\ref{imp9}) and (\ref{how}), \be
\la{go18}
 \Pr \{ L(K) \leq \se^\prime < U(K) \mid \se \} = \Pr
\{ L(K) < \se^{\prime \prime} \leq U(K) \mid \se \} > \f{\alpha +
\beta}{2}, \qqu \fa \se \in (\se^\prime, \se^{\prime \prime}). \ee
By virtue of (\ref{go18}) and recalling that both $\Pr \{ L(K) \leq
\se^\prime < U(K) \mid \se \}$ and $\Pr \{ L(K) < \se^{\prime
\prime} \leq U(K) \mid \se \}$ are continuous with respect to $\se
\in (\se^\prime, \se^{\prime \prime})$, we have
\[
\beta = \min \li ( \lim_{\se \downarrow \se^\prime}  \Pr \{ L(K)
\leq \se^\prime < U(K) \mid \se \}, \; \; \lim_{\se \uparrow
\se^{\prime \prime} } \Pr \{ L(K) < \se^{\prime \prime} \leq U(K)
\mid \se \} \ri ) \geq \f{\alpha + \beta}{2},
\]
leading to $\beta \geq \alpha$, which contradicts to $\alpha >
\beta$. Therefore, it must be true that $\alpha = \beta$.  That is,
{\small
\[ \inf_{\se \in [\se^\prime, \se^{\prime \prime}] } \Pr \{ L(K)
\leq \se \leq U(K) \mid \se \} =   \min \li (  \Pr \{ L(K) \leq
\se^\prime < U(K) \mid \se^\prime \}, \; \Pr \{ L(K) < \se^{\prime
\prime} \leq U(K) \mid \se^{\prime \prime}  \} \ri ). \]} It follows
that {\small \bel &  &  \inf_{\se \in [a, b]} \Pr \{ L (K) \leq \se
\leq
U (K) \mid \se \} \nonumber \\
&  &  = \min \{C(a), \; C_U(a), \; C(b), \; C_L(b) \} \cup \li \{
C_U (\se): \se \in \mscr{S}_U \cup \mscr{S}_L \ri \} \cup \li \{ C_L
(\se): \se \in \mscr{S}_L \cup \mscr{S}_U \ri \}. \la{cit} \eel}
Let
\[
\mscr{S}^\prime = \mscr{S}_U \cap \mscr{S}_L, \qqu \mscr{S}_U^\prime
= \mscr{S}_U \setminus \mscr{S}^\prime, \qqu \mscr{S}_L^\prime =
\mscr{S}_L \setminus \mscr{S}^\prime.
\]
Then,  \bel  &   & \li \{ C_U(\se): \se \in \mscr{S}_U \cup
\mscr{S}_L \ri \} \cup \li \{ C_L(\se): \se \in
\mscr{S}_U \cup \mscr{S}_L \ri \} \nonumber\\
 &  = & \li \{ C_U(\se): \se \in \mscr{S}_U^\prime \ri \} \cup \li \{ C_U(\se): \se \in \mscr{S}^\prime \ri \}
\cup \li \{ C_L(\se): \se \in \mscr{S}_L^\prime \ri \} \cup \li \{ C_L(\se): \se \in \mscr{S}^\prime \ri \} \nonumber\\
&   & \cup \li \{ C_U(\se): \se \in \mscr{S}_L^\prime \ri \} \cup
\li \{ C_L(\se): \se \in \mscr{S}^\prime_U \ri \}. \la{com1} \eel
For $\se \in \mscr{S}_U^\prime$, we have $ 0 \leq \Pr \{ L (K) = \se
< U (K) \mid \se \} \leq \Pr \{ L (K) = \se  \mid \se \} = 0$ and
thus \bee C_U(\se) - C_L(\se) & = & \Pr \{ L (K) \leq \se < U (K)
\mid \se \} - \Pr \{ L (K) < \se \leq U (K) \mid \se \}\\
& = & \Pr \{ L (K) = \se < U (K) \mid \se \} - \Pr \{ L (K) < \se = U (K) \mid \se \}\\
& = &  - \Pr \{ L (K) < \se = U (K) \mid \se \} \leq 0, \eee which
implies that \be \min \li \{ C_U(\se): \se \in \mscr{S}_U^\prime \ri
\} \leq \min \li \{ C_L(\se): \se \in \mscr{S}_U^\prime \ri \}.
\la{com2} \ee For $\se \in \mscr{S}_L^\prime$, we have $ 0 \leq \Pr
\{ L (K) < \se = U (K) \mid \se \} \leq \Pr \{ U (K) = \se \mid \se
\} = 0$ and thus \bee C_U(\se) - C_L(\se)
& = & \Pr \{ L (K) = \se < U (K) \mid \se \} - \Pr \{ L (K) < \se = U (K) \mid \se \}\\
& = & \Pr \{ L (K) = \se < U (K) \mid \se \} \geq 0, \eee which
implies that \be \min \li \{ C_U(\se): \se \in \mscr{S}_L^\prime \ri
\} \geq \min \li \{ C_L(\se): \se \in \mscr{S}_L^\prime \ri \}.
\la{com3} \ee Combing (\ref{com1}), (\ref{com2}) and (\ref{com3})
leads to \bel  &   & \min \li \{ C_U(\se): \se \in \mscr{S}_U \cup
\mscr{S}_L \ri \} \cup \li \{ C_L(\se): \se \in
\mscr{S}_U \cup \mscr{S}_L \ri \} \nonumber \\
 &  = & \min \li \{ C_U(\se): \se \in \mscr{S}_U^\prime \ri \} \cup \li \{ C_U(\se): \se \in \mscr{S}^\prime \ri \}
\cup \li \{ C_L(\se): \se \in \mscr{S}_L^\prime \ri \} \cup \li \{ C_L(\se): \se \in \mscr{S}^\prime \ri \} \nonumber\\
&  = & \min \li \{ C_U(\se): \se \in \mscr{S}_U \ri \} \cup \li \{
C_L(\se): \se \in \mscr{S}_L\ri \}, \la{cit2}
 \eel
 which implies that the minimum of the set $\{C(a), \; C_U(a), \; C(b), \; C_L(b) \} \cup  \{ C_U(\se): \se \in
\mscr{S}_U \cup \mscr{S}_L  \} \cup  \{ C_L(\se): \se \in \mscr{S}_U
\cup \mscr{S}_L  \}$ equals the minimum of $\{C(a), \; C_U(a), \;
C(b), \; C_L(b) \} \cup \li \{ C_U(\se): \se \in \mscr{S}_U \ri \}
\cup \li \{ C_L(\se): \se \in \mscr{S}_L \ri \}$.  This proves the
second statement of Theorem \ref{Fundamental}.

Clearly, the third statement of Theorem \ref{Fundamental} is already
justified in the course of proving the first two statements.  This
concludes the proof of Theorem \ref{Fundamental}.

\section{ Proof of Theorem \ref{infcov} } \la{infcov_app}

We shall first show that  $\inf_{p \in (0, 1)} \Pr \{ L(K) < p <
U(K) \mid p \}$ is equal to $\min_{p \in \mscr{S}} \Pr \{ L(K) < p <
U(K) \mid p \}$.  Clearly, as a consequence of  the assumption that
$0 = L(0) < U(0) < 1, \; 0 < L(n) < U(n) = 1$,  the sets
$\mscr{S}_L, \mscr{S}_U$ and $\mscr{S}$ are nonempty. Let $a$ and
$b$ be the minimum and maximum of $\mscr{S}$ respectively. Then, $0
< a \leq b < 1$ and $\inf_{p \in (0, 1)} \Pr \{ L(K) < p < U(K) \mid
p \}$ is equal to the minimum among $\inf_{p \in (0, a)} \Pr \{ L(K)
< p < U(K) \mid p \}, \; \inf_{p \in (b, 1)} \Pr \{ L(K) < p < U(K)
\mid p \}$ and $\inf_{p \in [a, b]} \Pr \{ L(K) < p < U(K) \mid p
\}$.  By the assumption that,  for any $\se \in (0, 1)$, there exist
two numbers $u$ and $v$ such that $\{ L(K) \leq \se \leq U(K) \} =
\{ u \leq K \leq v \}$, it can be easily shown by differentiation
that, for any $\vse \in (0, 1)$, $\Pr \{ L(K) \leq \vse \leq U(K)
\mid p \}$ is a continuous and unimodal function of $p \in (0, 1)$.
Hence, by virtue of Theorem \ref{Fundamental}, we have that \be
\la{need1}
 \inf_{p \in [a, b]}
\Pr \{ L(K) < p < U(K) \mid p \} = \min_{p \in \mscr{S}} \Pr \{ L(K)
< p < U(K) \mid p \}. \ee By Lemma 1 in Appendix A, we have that $\{
L(K) < p < U(K) \} = \{ L(K) < a \leq U(K) \} = \{ L(K) \leq
\f{a}{2} \leq U(K) \}$ for any $p \in (0, a)$.   By the assumption
that $0 = L(0) < U(0) < 1$, we have $U(0) \geq a > \f{a}{2}$, which
implies that $\{ K = 0 \} \subseteq \{ L(K) \leq \f{a}{2} \leq U(K)
\}$. Invoking the assumption that, for any $\se \in (0, 1)$, there
exist two numbers $u$ and $v$ such that $\{ L(K) \leq \se \leq U(K)
\} = \{ u \leq K \leq v \}$, we can conclude that there exists a
nonnegative integer $w$ such that $\{ L(K) \leq \f{a}{2} \leq U(K)
\} = \{ 0 \leq K \leq w \}$. Therefore, $\Pr \{ L(K) < p < U(K) \mid
p \} = \Pr  \{ L(K) < a \leq U(K) \mid p \} = \Pr \{ 0 \leq K \leq w
\mid p \}$ for any $p \in (0, a)$.  It can be easily shown that $\Pr
\{ 0 \leq K \leq w \mid p \}$ is monotonically decreasing with
respect to $p \in (0, a)$.  This implies that $\Pr \{ L(K) < a \leq
U(K) \mid p \} $ is monotonically decreasing with respect to $p \in
(0, a)$. Consequently, $\inf_{p \in (0, a)} \Pr \{ L(K) < p < U(K)
\mid p \} =  \inf_{p \in (0, a)} \Pr \{ L(K) < a \leq U(K) \mid p \}
= \lim_{p \uparrow a} \Pr \{ L(K) < a \leq U(K) \mid p \} = \Pr \{
L(K) < a \leq U(K) \mid a \}$ and it immediately follows that \be
\la{need2} \inf_{p \in (0, a)} \Pr \{ L(K) < p < U(K) \mid p \} \geq
\Pr \{ L(K) < a < U(K) \mid a \}. \ee By a similar argument, we can
show that \be \la{need3}
 \inf_{p \in (b,1)} \Pr \{ L(K) < p < U(K) \mid p
\} \geq \Pr \{ L(K) < b < U(K) \mid b \}. \ee Combining
(\ref{need1}), (\ref{need2}) and (\ref{need3}) leads to the
conclusion that $\inf_{p \in (0, 1)} \Pr \{ L(K) < p < U(K) \mid p
\}$ is equal to $\min_{p \in \mscr{S}} \Pr \{ L(K) < p < U(K) \mid p
\}$.

Next, we shall show that $\inf_{p \in (0, 1)} \Pr \{ L(K) \leq p
\leq U(K) \mid p \}$ is equal to the minimum of $\min_{p \in
\mscr{S}_L} \Pr \{ L(K) < p \leq U(K) \mid p \}$ and $\min_{p \in
\mscr{S}_U} \Pr \{ L(K) \leq p < U(K) \mid p \}$.  By a similar
argument as above, we can show that \be \la{need4} \inf_{p \in (0,
a)} \Pr \{ L(K) \leq p \leq U(K) \mid p \} = \Pr \{ L(K) < a \leq
U(K) \mid a \} < C \li ( \f{a}{2} \ri ) \ee and \be \la{need5}
\inf_{p \in (b, 1)} \Pr \{ L(K) \leq p \leq U(K) \mid p \} = \Pr \{
L(K) \leq b < U(K) \mid b \} < C \li ( \f{b + 1}{2} \ri ), \ee where
the notion of $C(.)$ is the same as that in Theorem
\ref{Fundamental}.

Let $Q_U$ denote the intersection of the interval $(\f{a}{2}, \f{b +
1}{2} )$ and the support of $U(K)$.  Let $Q_L$ denote the
intersection of the interval $(\f{a}{2}, \f{b + 1}{2} )$ and the
support of $L(K)$.  In the course of proving Theorem
\ref{Fundamental}, we have established (\ref{cit}). Invoking the
assumption that, for any $\se \in (0, 1)$, there exist two numbers
$u$ and $v$ such that $\{ L(K) \leq \se \leq U(K) \} = \{ u \leq K
\leq v \}$, we can conclude from (\ref{cit}) that $\inf_{p \in
[\f{a}{2}, \f{b + 1}{2} ] } \Pr \{ L(K) \leq p \leq U(K) \mid p \}$
is equal to the minimum of $\{ C( \f{a}{2} ), C ( \f{b + 1}{2} )
\}\cup \{ C_L (p): p \in Q_L \cup Q_U \} \cup \{ C_U (p): p \in Q_L
\cup Q_U \}$, where the meaning of $C(.), C_L(.), C_U(.) $ is the
same as that in Theorem \ref{Fundamental}.  Observing that
\[
\Pr \{ L(K) < a \leq U(K) \mid a \} \geq \min \{ C_L (p): p \in Q_L
\cup Q_U \}
\]
and
\[
\Pr \{ L(K) \leq b < U(K) \mid b \} \geq \min \{ C_U (p): p \in Q_L
\cup Q_U \},
\]
we have that the minimum among $\inf_{p \in (0, a)} \Pr \{ L(K) \leq
p \leq U(K) \mid p \}, \; \inf_{p \in (b, 1)} \Pr \{ L(K) \leq p
\leq U(K) \mid p \}$ and $\inf_{p \in [\f{a}{2}, \f{b + 1}{2} ] }
\Pr \{ L(K) \leq p \leq U(K) \mid p \}$ is equal to the minimum of
$\{ C_L (p): p \in Q_L \cup Q_U \} \cup \{ C_U (p): p \in Q_L \cup
Q_U \} = \{ C_L (p): p \in Q_L \} \cup \{ C_U (p): p \in  Q_U \} =
\{ C_L (p): p \in \mscr{S}_L \} \cup \{ C_U (p): p \in  \mscr{S}_U
\} $, where we have used (\ref{cit2}) established in the proof of
Theorem \ref{Fundamental}. It follows that the second statement of
Theorem \ref{infcov} on $\inf_{p \in (0, 1)} \Pr \{ L(K) \leq p \leq
U(K) \mid p \}$ holds true.

Finally, to show the third statement of Theorem \ref{infcov}, it is
sufficient to observe that $\Pr \{ L(K) < p = U(K) \mid p \} = 0$
for $p \in \mscr{S}_L$ and that $\Pr \{ L(K) = p < U(K) \mid p \} =
0$ for $p \in \mscr{S}_U$ as a consequence of  the assumption that
$\mscr{S_L} \cap \mscr{S}_L = \emptyset$.  The proof of Theorem
\ref{infcov} is thus completed.

\section{Proof of Theorem \ref{nonminima} } \la{nonminima_app}

For simplicity of notations, let $S_L = \{ L(k) \in (0, 1): 0 \leq k
\leq n \}$ and $S_U = \{ U(k) \in (0, 1): 0 \leq k \leq n  \}$. It
suffices to consider three exhaustive (but not mutually exclusive)
cases as follows.

Case (i):  $p \in S_L$;

Case (ii):  $p \in S_U$;

Case (iii):  $p \notin S_L \cup S_U$.

In Case (i), we can write $\{ L(K) \leq p \leq U(K) \} = \{ k \leq K
\leq l \}$, where $0 \leq k \leq l \leq n$ are integers.  Then,  $\{
L(K) \leq p - \ep \leq U(K) \} \subseteq \{ k \leq K \leq l - 1 \}$
for small enough $\ep > 0$. Thus, \bee &  & \Pr \{ L(K) \leq p  \leq
U(K) \mid p \} - \Pr \{ L(K) \leq
p - \ep \leq U(K) \mid p - \ep \}\\
&   &  = \Pr \{ k \leq K \leq l \mid p \} - \Pr \{ k \leq K \leq l
\mid p - \ep \} + \Pr \{  K = l \mid p - \ep \} \\
&  & \to \Pr \{  K = l \mid p \} > 0 \eee as $\ep \to 0$.   This
implies that $\Pr \{ L(K) \leq p  \leq U(K) \mid p \}$ is greater
than $\Pr \{ L(K) \leq p - \ep \leq U(K) \mid p - \ep \}$ for small
enough $\ep > 0$.  Hence, $p$ is not a local minima.

In Case (ii), we can write $\{ L(K) \leq p \leq U(K) \} = \{ k \leq
K \leq l \}$, where $0 \leq k \leq l \leq n$ are integers.  Then,
$\{ L(K) \leq p + \ep \leq U(K) \} \subseteq \{ k + 1 \leq K \leq l
\}$ for small enough $\ep
> 0$.  Thus,
\bee &  & \Pr \{ L(K) \leq p  \leq U(K) \mid p \} - \Pr \{ L(K) \leq
p + \ep \leq U(K) \mid p + \ep \}\\
&   &  = \Pr \{ k \leq K \leq l \mid p \} - \Pr \{ k \leq K \leq l
\mid p + \ep \} + \Pr \{  K = k \mid p + \ep \} \\
&  & \to \Pr \{  K = k \mid p \} > 0 \eee as $\ep \to 0$.   This
implies that $\Pr \{ L(K) \leq p  \leq U(K) \mid p \}$ is greater
than $\Pr \{ L(K) \leq p + \ep \leq U(K) \mid p + \ep \}$ for small
enough $\ep > 0$.  Hence, $p$ is not a local minima.

 In Case (iii), since $p \in (0,1)
\subseteq \cup_{k = 0}^{n} [L(k), U(k)]$, there must exist an
integer $k \in \{0, 1, \cd, n \}$ such that $p \in [ L(k), U(k)]$.
Thus, $\{ L(K) \leq p \leq U(K) \}$ is not an impossible event. As a
result, we can write $\{ L(K) \leq p \leq U(K) \} = \{ k \leq K \leq
l \}$, where $0 \leq k \leq l \leq n$ are integers.  Since $p \notin
S_L \cup S_U$,  we have that $\{ L(K) \leq p + \ep \leq U(K) \} = \{
k \leq K \leq l \}$ and $\{ L(K) \leq p - \ep \leq U(K) \} = \{ k
\leq K \leq l \}$ for small enough $\ep
> 0$.  Observing that $\Pr \{ k \leq K \leq l \mid \se \}$ is a
continuous and strictly monotone or unimodal function of $\se \in
(0, 1)$, we can conclude that $p$ is not a local minima.  The proof
of the theorem is thus completed.

\section{Proof of Theorem \ref{fundamental_thm9}} \la{secB}

For the simplicity of notations, define \[
\bi{m}{z} = \bec \f{ m! } { z! (m- z)! } & \tx{if} \; 0 \leq z \leq m,\\
0 & \tx{if} \; z < 0 \; \tx{or} \; z > m \eec
\]
for non-negative integer $m$ and arbitrary integer $z$.  We now
establish some preliminary results.

 \beL \label{tau
increase} Let $ 0 \leq M < N$. Define $T(k, M, N, n) = \li. {M
\choose k} {N-M-1 \choose n-k-1} \ri \slash {N \choose n}$. Then,
$\Pr \{ K \leq k \mid M \} - \Pr \{ K \leq k \mid M + 1 \} = T(k, M,
N, n)$ for any integer $k$. \eeL

\bpf

We first show the equation for $0 \leq k \leq M$.  We perform
induction on $k$. For $k = 0$, we have \bel \Pr \{ K \leq k \mid M
\} - \Pr \{ K \leq k \mid M + 1 \} & = & \Pr \{ K = 0 \mid M \} -
\Pr \{ K = 0 \mid M + 1 \} \nonumber\\
& = &   \frac{ {M \choose 0} {N-M \choose n} }
 { {N \choose n} } -  \frac{ {M + 1 \choose 0} {N-M-1 \choose n} }
 { {N \choose n} } \nonumber \\
 & =  & \frac{ {N-M-1 \choose n - 1} }
 { {N \choose n} } \la{comb}\\
 & =  & \frac{ {M \choose 0} {N-M-1 \choose n-0-1} }
 { {N \choose n} } = T(0, M, N, n) \nonumber, \eel where
 (\ref{comb}) follows from the fact that, for non-negative integer
 $m$, \be \la{basic} \bi{m + 1}{z+1} = \bi{m}{z} +
\bi{m}{z+1} \ee for any integer $z$.

Now suppose the lemma is true for $k-1$ with $1 \leq k \leq M$, i.e.,
 \[
\Pr \{ K \leq k - 1 \mid M \} - \Pr \{ K \leq k - 1 \mid M + 1 \} =
\frac{ {M \choose k-1} {N-M-1 \choose n-k} }
 { {N \choose n} }.
 \]
Then, {\small \bel  \Pr \{ K \leq k \mid M \} - \Pr \{ K \leq k \mid
M + 1 \}
 & = &  \Pr \{ K \leq k - 1 \mid M \} - \Pr \{ K \leq k - 1 \mid M + 1 \} \nonumber\\
 &   &  + \frac{ {M \choose k} {N-M
\choose n - k} } { {N \choose n} } -  \frac{ {M + 1 \choose k} {N-M-1 \choose n - k} }
 { {N \choose n} } \nonumber\\
 & = & \frac{ {M \choose k-1} {N-M-1 \choose n-k} }
 { {N \choose n} } +  \frac{ {M \choose k} {N-M
\choose n - k} } { {N \choose n} } -  \frac{ {M + 1 \choose k} {N-M-1 \choose n - k} }
 { {N \choose n} } \nonumber\\
 &  =  & \frac{ {M \choose k} {N-M
\choose n - k} } { {N \choose n} } - \li [ \frac{ {M + 1 \choose k} {N-M-1 \choose n - k} }
 { {N \choose n} } - \frac{ {M \choose k-1} {N-M-1 \choose n-k} }
 { {N \choose n} } \ri ] \nonumber\\
& = & \frac{ {M \choose k} {N-M \choose n - k} } { {N \choose n} } - \frac{ {M \choose k} {N-M - 1 \choose n -
k} } { {N \choose n} } \la{eqaa}\\
& = & \frac{ {M \choose k} {N-M - 1 \choose n - k - 1} } { {N \choose n} } \la{eqbb} \eel} where (\ref{eqaa})
and (\ref{eqbb}) follows from (\ref{basic}).  Therefore, we have shown the lemma for $0 \leq k \leq M$.

For $k > M$, we have $\Pr \{ K \leq k \mid M \} = \Pr \{ K \leq k
\mid M + 1 \}  = 1$ and $T(k, M, N, n) = 0$.  For $k < 0$, we have
$\Pr \{ K \leq k \mid M \} = \Pr \{ K \leq k \mid M + 1 \} = 0$ and
$T(k, M, N, n) = 0$. Thus, the lemma is true for any integer $k$.

\epf

\beL \la{dou} Let $1 \leq M \leq N$ and $k \leq l$. Then,
\[ \Pr \{ k \leq K \leq l \mid M \}  - \Pr \{ k \leq K \leq l \mid M - 1 \} = T(k - 1, M-1, N, n) - T(l, M - 1, N, n).
\]
\eeL

\bpf

To show the lemma, it suffices to consider $6$ cases as follows.

Case (i): $0 < n < k \leq l$. In this case, $\Pr \{ k \leq K \leq l
\mid M \}  = \Pr \{ k \leq K \leq l \mid M - 1 \} = 0$ and $T(k - 1,
M-1, N, n) = T(l, M - 1, N, n) = 0$.

Case (ii): $k \leq l < 0 < n$. In this case, $\Pr \{ k \leq K \leq l
\mid M \}  = \Pr \{ k \leq K \leq l \mid M - 1 \}  = 0$ and $T(k -
1, M-1, N, n) = T(l, M - 1, N, n) = 0$.

Case (iii): $k \leq 0 < n \leq l$.  In this case, $\Pr \{ k \leq K
\leq l \mid M \}  = \Pr \{ k \leq K \leq l \mid M - 1 \} = 1$ and
$T(k - 1, M-1, N, n) = T(l, M - 1, N, n) = 0$.

Case (iv): $k \leq 0 \leq l < n$.  In this case, $T(k - 1, M-1, N,
n) = 0$ and, by Lemma \ref{tau increase}, {\small \bee \Pr \{ k \leq
K \leq l \mid M \}  - \Pr \{ k \leq K \leq l \mid M - 1 \} & = & \Pr
\{ K \leq l \mid M \}  - \Pr \{ K \leq l \mid M - 1 \}\\
& = & T(k - 1, M-1, N, n) - T(l, M - 1, N, n).
 \eee}

Case (v): $0 < k \leq n \leq l$.  In this case, $T(l, M-1, N, n) =
0$ and, by Lemma \ref{tau increase}, {\small \bee \Pr \{ k \leq K
\leq l \mid M \}  - \Pr \{ k \leq K \leq l \mid M - 1 \}
& = & \Pr \{ K < k \mid M - 1 \}  - \Pr \{ K < k \mid M \}\\
& = & T(k - 1, M-1, N, n) - T(l, M - 1, N, n).
 \eee}

Case (vi): $0 < k \leq l < n$.  In this case, by Lemma \ref{tau
increase},  {\small \bee &  & \Pr \{ k \leq K \leq l \mid M \}  -
\Pr \{
k \leq K \leq l \mid M - 1 \}\\
 & = & [ \Pr \{ K \leq l \mid M \} - \Pr \{ K < k \mid M \}] - [\Pr \{ K \leq l \mid M - 1 \} - \Pr \{ K < k \mid M - 1 \} ]\\
& = & [ \Pr \{ K \leq l \mid M \}  - \Pr \{ K \leq l \mid M - 1 \} ]  - [ \Pr \{ K < k \mid M \} - \Pr \{ K < k \mid M - 1 \} ]\\
& = & T(k - 1, M-1, N, n) - T(l, M - 1, N, n).
 \eee}

\epf

\beL \la{lem22} Let $l \geq 0$ and $k < n$.  Then, $\li \lf \f{ n M} { N+1} \ri \rf \geq l$ for $M \geq 1 + \li
\lf \f{ N l } { n - 1 } \ri \rf$, and  $\li \lf \f{  n M} { N+1} \ri \rf \leq k - 1$ for $M \leq 1 + \li \lf \f{
N (k - 1) } { n - 1 } \ri \rf$.

\eeL

\bpf

 To show the first part of the lemma, observe that $(N+1 - n) l \geq 0$,
 by which we can show $\f{n N l}{n-1}
\geq (N+1) l$. Hence, $n \li ( 1 + \li \lf \f{ N l } { n - 1 } \ri \rf \ri )
> \f{n N l}{n-1} \geq (N+1)
l$.  That is, $ \f{n} { N+1} \li ( 1 + \li \lf \f{ N l } { n - 1 } \ri \rf \ri ) >  l$.  It follows that $\li
\lf \f{n}{ N+1} \li ( 1 + \li \lf \f{ N l } { n - 1 } \ri \rf \ri ) \ri \rf \geq l$. Since the floor function is
non-decreasing, we have $\li \lf \f{ n M} { N+1} \ri \rf \geq l$ for $M \geq 1 + \li \lf \f{ N l } { n - 1 } \ri
\rf$.

To prove the second part of the lemma, note that $(N + 1 - n) (n - k) > 0$, from which we can deduce $1 + \f{
N(k - 1) } { n - 1 } < \f{ (N +1) k } { n }$.  Hence, $1 + \li \lf \f{ N(k - 1) } { n - 1 } \ri \rf < \f{ (N +1)
k  } { n }$, i.e., $\f{ n } {N + 1} \li ( 1 + \li \lf \f{ N(k - 1) } { n - 1 } \ri \rf \ri ) < k$, leading to
$\li \lf \f{ n } {N + 1} \li ( 1 + \li \lf \f{ N(k - 1) } { n - 1 } \ri \rf \ri ) \ri \rf \leq k - 1$. Since the
floor function is non-decreasing, we have $\li \lf \f{ n M} { N+1} \ri \rf \leq k - 1$ for $M \leq 1 + \li \lf
\f{ N (k - 1) } { n - 1 } \ri \rf$.

\epf

\beL \la{monotone}
 Let $0 \leq r \leq n$.  Then, the following statements hold true.

(I) \[ T(r - 1, M - 1, N, n) \leq T(r, M - 1, N, n) \qu \mrm{for} \qu 1 \leq r \leq \li \lf \f{ n M} { N+1} \ri
\rf;
\]
\[ T(r + 1, M - 1, N, n) \leq T(r, M - 1, N, n) \qu \mrm{for}
\qu  \li \lf \f{ n M} { N+1} \ri \rf \leq r \leq n - 1.
\]

(II)
\[
T(r, M - 2, N, n) \leq T(r, M - 1, N, n) \qu \mrm{for} \qu 1 < M \leq 1 + \li \lf \f{ N r } { n - 1 } \ri \rf;
\]
\[
T(r, M, N, n) \leq T(r, M - 1, N, n) \qu \mrm{for} \qu  1 + \li \lf \f{ N r } { n - 1 } \ri \rf \leq M < N.
\]

\eeL

\bpf

To show statement (I), note that $T(r, M -1, N, n) = 0$ for $\min(M - 1, n-1) < r \leq n$.  Our calculation
shows that
 \[
 \f{ T(r - 1, M -1, N, n) } { T(r, M -1, N, n)
} = \f{r}{ M - r} \f{ N - M + 1 - (n -r) } { n - r } \leq 1 \qu \mrm{for} \qu 1 \leq  r \leq \f{ n M } { N + 1 }
\]
 and
 \[
 \f{ T(r - 1, M-1, N, n) }
{ T(r, M -1, N, n) } > 1 \qu \mrm{for} \qu \f{ n M } { N + 1 } < r \leq \min(M - 1, n-1).
\]

To show statement (II), note that $T(r, M - 1, N, n) = 0$ for $1 \leq M < r + 1$, and $T(r, M - 1, N, n) \geq
T(r, M - 2, N, n) = 0$ for $M = r + 1$.  Direct computation shows that {\small \[ \f{ T(r, M -1, N, n) } { T(r,
M - 2, N, n) } = \f{ M -1 } {M - 1 - r}
 \f{ N - M  + 2 - (n -r) } {N - M + 1} \geq 1 \qu \mrm{ for} \qu r + 1 < M  \leq 1 +  \f{ N r } { n - 1
 },
 \]}
  and
  \[
  \f{ T(r, M - 1, N, n) } { T(r, M - 2, N,
n) } < 1 \qu \mrm{ for } \qu  1 + \f{ N r } { n - 1 } < M \leq N.
\]

\epf

\beL \la{lem77} Let $0 \leq \mcal{L} \leq \mcal{U} \leq N$.  Then,
for any integers $k$ and $l$,  $\Pr \{ k \leq K \leq l \mid M \}$ is
unimodal with respect to $M$ for $\mcal{L} \leq M \leq \mcal{U}$.
\eeL

\bpf

Clearly, the lemma is trivially true if $k > l$.  Hence, to show the lemma, it suffices to consider $6$ cases as
follows.

Case (i): $0 < n < k \leq l$. In this case, $\Pr \{ k \leq K \leq l
\mid M \} = 0$ for any $M \in [\mcal{L}, \mcal{U}]$.

Case (ii): $k \leq l < 0 < n$. In this case, $\Pr \{ k \leq K \leq l
\mid M \} = 0$ for any $M \in [\mcal{L}, \mcal{U}]$.

Case (iii): $k \leq 0 < n \leq l$. In this case, $\Pr \{ k \leq K
\leq l \mid M \} = 1$ for any $M \in [\mcal{L}, \mcal{U}]$.

Case (iv): $k \leq 0 \leq l < n$.  In this case, $\Pr \{ k \leq K
\leq l \mid M \} = \Pr \{ K \leq l \mid M \}$ is non-increasing with
respect to $M \in [\mcal{L}, \mcal{U}]$ as can be seen from Lemma
\ref{tau increase}.

Case (v): $0 < k \leq n \leq l$. In this case, $\Pr \{ k \leq K \leq
l \mid M \}  = 1 - \Pr \{ K < k \mid M \}$ is non-decreasing with
respect to $M \in [\mcal{L}, \mcal{U}]$ as can be seen from Lemma
\ref{tau increase}.

Clearly, the lemma is true for the above five cases.

Case (vi): $0 < k \leq l < n$.  Define {\small $\De(k, l, M, N, n) =
\Pr \{ k \leq K \leq l \mid M \}  - \Pr \{ k \leq K \leq l \mid M -
1 \}$}. By Lemma \ref{dou}, $\De(k, l, M, N, n) = T(k - 1, M-1, N,
n) - T(l, M - 1, N, n)$.

Invoking Lemma \ref{lem22}, for $M \geq 1 + \li \lf \f{ N l } { n - 1 } \ri \rf$, we have that $\li \lf \f{  n
M} { N+1} \ri \rf \geq l$ and thus, by statement (I) of Lemma \ref{monotone}, $T(r, M - 1, N, n)$ is
non-decreasing with respect to $r \leq l$. Consequently, $T(k - 1, M -1, N, n) \leq T(l, M -1, N, n)$, leading
to $\De(k, l, M, N, n) \leq 0$ for $M \geq 1 + \li \lf \f{ N l } { n - 1 } \ri \rf$.

Similarly, applying Lemma \ref{lem22}, for $M \leq 1 + \li \lf \f{ N (k - 1) } { n - 1 } \ri \rf$, we have that
$\li \lf \f{  n M} { N+1} \ri \rf \leq k - 1$ and thus, by statement (I) of Lemma \ref{monotone}, $T(r, M - 1,
N, n)$ is non-increasing with respect to $r \geq k - 1$. Consequently, $T(k -1, M -1, N, n) \geq T(l, M -1, N,
n)$, leading to $\De(k, l, M, N, n) \geq 0$ for $M \leq 1 + \li \lf \f{ N (k - 1) } { n - 1 } \ri \rf$.

By statement (II) of Lemma \ref{monotone}, for $1 + \li \lf \f{ N (k
- 1) } { n - 1 } \ri \rf \leq M \leq 1 + \li \lf \f{ N l } { n - 1 }
\ri \rf$, we have that $T(l, M - 1, N, n)$ is non-decreasing with
respect to $M$ and that $T(k - 1, M - 1, N, n)$ is non-increasing
with respect to $M$. It follows that $\De(k, l, M, N, n)$ is
non-increasing with respect to $M$ in this range.  Therefore, there
exists an integer $M^*$ such that $1 + \li \lf \f{ N (k - 1) } { n -
1 } \ri \rf \leq M^* \leq 1 + \li \lf \f{ N l } { n - 1 } \ri \rf$
and that $\De(k, l, M, N, n) \geq 0$ for $0 \leq M \leq M^*$, and
$\De(k, l, M, N, n) \leq 0$ for $M^* \leq M \leq N$.  This implies
that $\Pr \{ k \leq K \leq l \mid M \}$ is non-decreasing for $0
\leq M \leq M^*$ and non-increasing for $M^* \leq M \leq N$. This
concludes the proof of the lemma.

\epf

\beL \la{compa} Let $0 \leq M <  N$. Then,  $\Pr \{ g \leq K \leq h
+ 1 \mid M + 1 \} \geq \Pr \{ g \leq K \leq h \mid M \}$ for any
integers $g$ and $h$. \eeL

\bpf Clearly, the lemma is trivially true if $g > h$.  Hence, to
show the lemma, it suffices to consider the case $g \leq h$.  Note
that, by Lemma \ref{dou}, \bee &   & \Pr \{ g \leq K \leq h + 1 \mid
M + 1 \} - \Pr \{ g \leq K \leq h \mid M \}\\
 & = & \bi{M +1}{h+1} \bi{N - M - 1}{n - h - 1} \li \sh \bi{N}{n} \ri.
 + \Pr \{ g \leq K \leq h \mid M + 1 \} - \Pr \{ g \leq K \leq h \mid M \}\\
& = &  \bi{M +1}{h+1} \bi{N - M - 1}{n - h - 1}
\li \sh \bi{N}{n} \ri. + T(g - 1, M, N, n) - T(h, M, N, n)\\
& = & \li [ \bi{M +1}{h+1} \bi{N - M - 1}{n - h - 1}  -
\bi{M}{h} \bi{N - M - 1}{n - h - 1} \ri ] \li \sh \bi{N}{n} \ri. + T(g - 1, M, N, n)\\
& = &  \bi{M}{h+1} \bi{N - M - 1}{n - h - 1} \li \sh \bi{N}{n} \ri.
+ T(g - 1, M, N, n) \geq  0, \eee where the last equality follows
from (\ref{basic}).

 \epf

\beL \la{compb} Let $0 < M \leq N$. Then, $\Pr \{ g - 1 \leq K \leq
h \mid M - 1 \} \geq \Pr \{ g \leq K \leq h \mid M \}$ for any
integers $g$ and $h$. \eeL

\bpf

Clearly, the lemma is trivially true if $g > h$.  Hence, to show the
lemma, it suffices to consider the case $g \leq h$.  Note that, by
Lemma \ref{dou}, \bee &   & \Pr \{ g - 1 \leq K
\leq h \mid M - 1 \} - \Pr \{ g \leq K \leq h \mid M \}\\
& = & \bi{ M - 1 } {g - 1  } \bi{ N - M + 1  } { n - g + 1 } \li \sh
\bi{N}{n} \ri. + \Pr \{ g \leq K
\leq h \mid M - 1 \} - \Pr \{ g \leq K \leq h \mid M \}\\
& = &  \bi{ M - 1 } {g - 1  } \bi{ N - M + 1  } { n - g + 1
} \li \sh \bi{N}{n} \ri. + T(h, M - 1, N, n) - T(g - 1, M - 1, N, n)\\
& = &  \li [ \bi{ M - 1 } {g - 1  } \bi{ N - M + 1  } { n - g + 1 }
- \bi{ M - 1 } {g - 1  } \bi{ N
- M } { n - g } \ri ] \li \sh \bi{N}{n} \ri. + T(h, M - 1, N, n)\\
& = & \bi{ M - 1 } {g - 1  } \bi{ N - M } { n - g + 1 } \li \sh
\bi{N}{n} \ri. + T(h, M - 1, N, n) \geq 0,
 \eee
 where the last equality
follows from (\ref{basic}).

\epf

\beL \la{last}  Suppose that $\{ M^\prime < L(K) <  M^{\prime
\prime} \} = \{ M^\prime < U(K) < M^{\prime \prime} \}  =
\emptyset$.  Then, $\Pr \{ L(K) < M < U(K) \mid M \}$ is unimodal
with respect to $M$ for $M^\prime \leq M \leq M^{\prime \prime}$.
\eeL

\bpf

First, we shall show the following facts:

(i) If $\{ L(K) = M^\prime \} = \emptyset$, then $\{ L(K) < M \}  =
\{ L(K) < M^\prime \} = \{ L(K) < M^{\prime \prime} \}$ for
$M^\prime \leq M \leq M^{\prime \prime}$.

(ii) If $\{ L(K) = M^\prime \} \neq \emptyset$, then $\{ L(K) < M \}
= \{ L(K) \leq M^\prime \} = \{ L(K) < M^{\prime \prime} \}$ for
$M^\prime < M \leq M^{\prime \prime}$.

(iii) If $\{ U(K) = M^{\prime \prime} \} = \emptyset$, then $\{ U(K)
> M \} = \{ U(K) > M^\prime \} = \{ U(K) > M^{\prime \prime} \}$
for $M^\prime \leq M \leq M^{\prime \prime}$.

(iv) If $\{ U(K) = M^{\prime \prime} \} \neq \emptyset$, then $\{
U(K) > M \} = \{ U(K) > M^\prime \} = \{ U(K) \geq M^{\prime \prime}
\}$ for $M^\prime \leq M < M^{\prime \prime}$.

\bsk

To show statement (i), making use of $\{ L(K) = M^\prime \} = \{
M^\prime < L(K) < M^{\prime \prime} \} = \emptyset$, we have $\{
M^\prime \leq L(K) < M \} = \{ M^\prime < L(K) < M \} \subseteq \{
M^\prime < L(K) < M^{\prime \prime} \} = \emptyset$ and $\{ L(K) < M
\} = \{ L(K) < M^\prime \} \cup \{ M^\prime  \leq L(K) < M \} = \{
L(K) < M^\prime \}$ for $M^\prime \leq M \leq M^{\prime \prime}$. On
the other hand, $\{ L(K) < M \} = \{ L(K) < M^{\prime \prime} \}
\setminus \{ M \leq L(K) < M^{\prime \prime} \} = \{ L(K) <
M^{\prime \prime} \}$ for $M^\prime \leq M \leq M^{\prime \prime}$.

\bsk

To show statement (ii), making use of $\{ M^\prime < L(K) <
M^{\prime \prime} \} = \emptyset$, we have $\{ M^\prime < L(K) < M
\} \subseteq \{ M^\prime < L(K) < M^{\prime \prime} \} = \emptyset$
and $\{ L(K) < M \} = \{ L(K) \leq M^\prime \} \cup \{ M^\prime  <
L(K) < M \} = \{ L(K) \leq M^\prime \}$ for $M^\prime \leq M \leq
M^{\prime \prime}$. On the other hand, $\{ L(K) < M \} = \{ L(K) <
M^{\prime \prime} \} \setminus \{ M \leq L(K) < M^{\prime \prime} \}
= \{ L(K) < M^{\prime \prime} \}$ for $M^\prime < M \leq M^{\prime
\prime}$.

\bsk

To show statement (iii), using $\{ U(K) = M^{\prime \prime} \} = \{
M^\prime < U(K) < M^{\prime \prime} \} = \emptyset$,  we have $\{
M^\prime < U(K) \leq M  \} \subseteq \{ M^\prime < U(K) \leq
M^{\prime \prime} \} = \emptyset$ and $\{ U(K) > M \} = \{ U(K) >
M^\prime \} \setminus \{  M^\prime < U(K) \leq M  \} = \{ U(K) >
M^\prime \}$ for $M^\prime \leq M \leq M^{\prime \prime}$.  On the
other hand, $\{ U(K) > M \} = \{ U(K) > M^{\prime \prime} \} \cup \{
M < U(K) \leq M^{\prime \prime} \} = \{ U(K) > M^{\prime \prime} \}$
for $M^\prime \leq M \leq M^{\prime \prime}$.

\bsk

To show statement (iv), note that $\{ U(K) > M \} = \{ U(K) >
M^\prime \}$ for $M^\prime \leq M < M^{\prime \prime}$.  On the
other hand, $\{ U(K) > M \} = \{ U(K) \geq M^{\prime \prime} \} \cup
\{ M < U(K) < M^{\prime \prime} \} = \{ U(K) \geq  M^{\prime \prime}
\}$ for $M^\prime \leq M < M^{\prime \prime}$.

\bsk

Now, to show the lemma, it suffices to consider four cases as
follows.

\bed

\item Case (i): $\{ L(K) = M^\prime \} = \emptyset, \; \{ U(K) = M^{\prime \prime} \} = \emptyset$.

\item Case (ii): $\{ L(K) = M^\prime \} = \emptyset, \; \{ U(K) = M^{\prime \prime} \} \neq \emptyset$.

\item Case (iii): $\{ L(K) = M^\prime \} \neq \emptyset, \; \{ U(K) = M^{\prime \prime} \} = \emptyset$.

\item Case (iv): $\{ L(K) = M^\prime \} \neq \emptyset, \; \{ U(K) = M^{\prime \prime} \} \neq \emptyset$.

\eed

In Case (i), making use of facts (i) and (iii), we have $\{ L(K) < M
< U(K) \} = \{ L(K) < M^\prime < U(K) \} = \{ L(K) < M^{\prime
\prime} < U(K) \}$ for $M^\prime \leq M \leq M^{\prime \prime}$.
Invoking Lemma \ref{lem77}, we have that  $\Pr \{ L(K) < M < U(K)
\mid M \}$ is unimodal with respect to $M$ for $M^\prime \leq M \leq
M^{\prime \prime}$.

\bsk

In Case (ii), making use of facts (i) and (iv), we have  $\{ L(K) <
M < U(K) \} = \{ L(K) < M^\prime < U(K) \} = \{ L(K) < M^{\prime
\prime} \leq U(K) \}$ for $M^\prime \leq M < M^{\prime \prime}$.
Invoking Lemma \ref{lem77},  we have that $\Pr \{ L(K) < M < U(K)
\mid M \}$ is unimodal with respect to $M$ for $M^\prime \leq M <
M^{\prime \prime}$.  Since $\{ M^{\prime \prime} = U(K) \} \neq
\emptyset$ and $U(K)$ is monotonically increasing,  we have $\{
M^{\prime \prime} \leq U(K) \} = \{ K \geq \udl{k} \}$ and $\{
M^{\prime \prime} < U(K) \} = \{ K \geq \ovl{k} + 1 \}$,  where
$\udl{k} = \min \{ k : U(k) \geq M^{\prime \prime} \} \leq \ovl{k} =
\max \{ k : U(k) \leq M^{\prime \prime} \}$.  Therefore, as a result
of Lemma \ref{compb}, \be \la{use1}
 \Pr \{ L(K) <
M^{\prime \prime} \leq U(K) \mid M^{\prime \prime} - 1 \} \geq \Pr
\{ L(K) < M^{\prime \prime} < U(K) \mid M^{\prime \prime} \}. \ee It
follows that $\Pr \{ L(K) < M < U(K) \mid M \}$ is unimodal with
respect to $M$ for $M^\prime \leq M \leq M^{\prime \prime}$.

In Case (iii), making use of facts (ii) and (iii), we have $\{ L(K)
< M < U(K) \} = \{ L(K) \leq M^\prime < U(K) \} = \{ L(K) <
M^{\prime \prime} < U(K) \}$ for $M^\prime < M \leq M^{\prime
\prime}$. Invoking Lemma \ref{lem77}, we have that $\Pr \{ L(K) < M
< U(K) \mid M \}$ is unimodal with respect to $M$ for $M^\prime < M
\leq M^{\prime \prime}$.  Since $\{ M^\prime = L(K) \} \neq
\emptyset$ and $L(K)$ is monotonically increasing,  we have $\{
M^\prime \geq L(K) \} = \{ K \leq \ovl{k} \}$ and $\{ M^\prime >
L(K) \} = \{ K \leq \udl{k} - 1 \}$,  where $\udl{k} = \min \{ k :
L(k) \geq M^\prime \} \leq \ovl{k} = \max \{ k : L(k) \leq M^\prime
\}$.  Therefore, as a result of Lemma \ref{compb}, \be \la{use2} \Pr
\{ L(K) < M^\prime < U(K) \mid M^\prime \} \leq \Pr \{ L(K) \leq
M^\prime < U(K) \mid M^\prime + 1 \}. \ee It follows that $\Pr \{
L(K) < M < U(K) \mid M \}$ is unimodal with respect to $M$ for
$M^\prime \leq M \leq M^{\prime \prime}$.

In Case (iv), making use of facts (ii) and (iv), we have $\{ L(K) <
M < U(K) \} = \{ L(K) \leq M^\prime < U(K) \} = \{ L(K) < M^{\prime
\prime} \leq U(K) \}$ for $M^\prime < M < M^{\prime \prime}$.
Invoking Lemma \ref{lem77}, we have that $\Pr \{ L(K) < M < U(K)
\mid M \}$ is unimodal with respect to $M$ for $M^\prime < M <
M^{\prime \prime}$.  Recalling (\ref{use1}) and (\ref{use2}), we
have that  $\Pr \{ L(K) < M < U(K) \mid M \}$ is unimodal with
respect to $M$ for $M^\prime \leq M \leq M^{\prime \prime}$.

\epf

Finally,  we are in a position to prove the Theorem
\ref{fundamental_thm9}.  Let $M^\prime < M^{\prime \prime}$ be two
consecutive distinct elements of $I_{UL}$. Then, $\{ M^\prime < L(K)
< M^{\prime \prime} \} = \{ M^\prime < U(K) < M^{\prime \prime} \} =
\emptyset$.  By Lemma \ref{last}, we have that $\Pr \{ L(K) < M <
U(K) \mid M \}$ is unimodal with respect to $M$ for $M^\prime \leq M
\leq M^{\prime \prime}$. Since this argument holds for any
consecutive distinct elements of the set $I_{UL}$, Theorem
\ref{fundamental_thm9} is established.

\end{document}